\documentclass{amsart}
\usepackage{amsmath,amssymb}
\usepackage[all]{xy}

\theoremstyle{plain}
\newtheorem{theorem}{Theorem}[section]
\newtheorem{proposition}[theorem]{Proposition}
\newtheorem{lemma}[theorem]{Lemma}
\newtheorem{corollary}[theorem]{Corollary}

\newtheorem{question}[theorem]{Question}

\theoremstyle{definition}
\newtheorem{definition}[theorem]{Definition}
\newtheorem{example}[theorem]{Example}
\newtheorem{examples}[theorem]{Examples}
\newtheorem{notation}[theorem]{Notation}

\theoremstyle{remark}

\newtheorem{remarks}[theorem]{Remarks}

\newcommand{\secref}[1]{Section~\ref{#1}}
\newcommand{\thmref}[1]{Theorem~\ref{#1}}
\newcommand{\propref}[1]{Proposition~\ref{#1}}
\newcommand{\lemref}[1]{Lemma~\ref{#1}}
\newcommand{\corref}[1]{Corollary~\ref{#1}}

\newcommand{\remsref}[1]{Remarks~\ref{#1}}
\newcommand{\defref}[1]{Definition~\ref{#1}}
\newcommand{\exsref}[1]{Examples~\ref{#1}}

\def\R{{\mathbb \R}}

\def\Q{{\mathbb Q}}

\def\map{\mathrm{map}}

\def\cat0{\mathrm{cat}_0}
\def\cat{\mathrm{cat}}

\def\catBG{\mathrm{cat}_\mathrm{BG}}
\def\catAL{\mathrm{cat}}
\def\catFH{\mathrm{cat}_\mathrm{FH}}
\def\cocatBG{\mathrm{cocat}_\mathrm{BG}}
\def\cocatAL{\mathrm{cocat}}
\def\id{\mathrm{id}}
\def\triv{*}

\begin{document}

\title[Homotopy Actions, Cyclic Maps and their Duals]{Homotopy Actions, Cyclic Maps and their Duals}

\author{Martin Arkowitz}

\address{Department of Mathematics,
          Dartmouth College,
          Hanover NH 03755}

\email{M.Arkowitz@Dartmouth.edu}

\author{Gregory  Lupton}

\address{Department of Mathematics,
          Cleveland State University,
          Cleveland OH 44115}

\email{G.Lupton@csuohio.edu}

\date{\today}

\keywords{Action, cyclic map, category of a
map, coaction, cocyclic map, cocategory of a
map}

\subjclass[2000]{55Q05, 55M30, 55P30}

\begin{abstract}

An \emph{action of $A$ on $X$} is a map $F\colon A\times X \to X$
such that $F\vert_X = \id \colon X\to X$.  The restriction $F\vert_A
\colon A \to X$ of an action is called a \emph{cyclic map}. Special
cases of these notions include group actions and the Gottlieb groups
of a space, each of which has been studied extensively. We prove
some general results about actions and their Eckmann-Hilton duals.
For instance, we classify the actions on an $H$-space that are
compatible with the $H$-structure. As a corollary, we prove that if
any two actions $F$ and $F'$ of $A$ on $X$ have cyclic maps $f$ and
$f'$ with $\Omega f = \Omega f'$, then $\Omega F$ and $\Omega F'$
give the same action of $\Omega A$ on $\Omega X$. We introduce a new
notion of the category of a map $g$ and prove that $g$ is cocyclic
if and only if the category is less than or equal to $1$.  From this
we conclude that if $g$ is cocyclic, then the Berstein-Ganea
category of $g$ is $\le 1$.  We also briefly discuss the
relationship between a map being cyclic and its cocategory being
$\le 1$.
\end{abstract}

\maketitle

\section{Introduction}\label{sec:intro}

Group actions on spaces are natural objects of study in topology.
Their analogues in homotopy theory lead to the notion of a homotopy
action of one space on another.  The restriction to the basepoint of
$X$ of a group action $G \times X \to X$ yields a map $G \to X$,
known as the \emph{orbit map} of the action.  These maps, and their
homotopy-theoretic analogues, are frequently of interest. For
instance, the \emph{$n$th Gottlieb group} of a space $X$, denoted
$G_n(X)$, may be defined as homotopy classes of maps $f \colon S^n
\to X$ such that $(f \mid \id) \colon S^n \vee X \to X$ admits an
extension $ S^n\times X \to X$ up to homotopy \cite{Go1}. That is,
the $n$th Gottlieb group is the set of homotopy classes of orbit
maps of $S^n$-actions on $X$. These groups can be
generalized.  For example, in \cite{Var} the homotopy set $G(A,
X)$ of \emph{cyclic maps} from $A$ to $X$ is defined as the homotopy
classes of maps $f \colon A \to X$ such that $(f \mid \id) \colon A
\vee X \to X$ admits an extension $ A \times X \to X$.  Such an extension is called an $f$-action, and so a cyclic map $f:A \to X$
is just the orbit map of an $f$-action of $A$ on $X$. Actions and cyclic
maps have been well-studied in the literature.  For the former, see
\cite{Gan67, Zab76, Iw-Od}. For the latter, see \cite{Var, 
Hal-Var, Lee-Woo}. The Gottlieb groups have also received much
attention in rational homotopy, following the results of
F{\'e}lix-Halperin in \cite{F-H}.

In this paper we prove some general results about actions and cyclic
maps.  In \secref{sec:actions}, we focus on actions and the
operations on homotopy sets that are induced by actions.  We deduce
several formulas, such as one for the homomorphism induced on
homotopy sets (\lemref{lem:cyclic on homotopy}).  We also conclude
that any two actions loop to the same action (\propref{prop:loop of
an action}), which implies that any two multiplications on
an $H$-space loop to the same multiplication (\corref{cor:loop m =
loop m'}). Our results dualize, and we summarize the dual results at
the end of \secref{sec:actions}. In \secref{sec:category}, we concentrate
on the dual of a cyclic map, called a \emph{cocyclic map}.  We
introduce a new version of the Lusternik-Schnirelmann category
of a map and show that the maps of category $\le 1$ are precisely
the cocyclic ones (Theorem 3.7).
Since there has been increasing interest in Lusternik-Schnirelmann category recently (see \cite{CLOT}), this gives motivation for
the study of cocyclic maps.  Furthermore, we clarify the relation
between cocyclic maps and connecting maps of cofibration sequences (Theorem 
3.16).
In doing so, we fill a gap in the literature.  The relation between
cyclic maps and connecting maps of fibration sequences has long been
well-understood \cite{Go1, Hal-Var}.  The dual situation, however,
has remained unclear due to the lack of a dual for the classifying
space construction.

All spaces, maps and homotopies are assumed to be based. All spaces
are assumed to be CW complexes. We use $\triv\colon X\to Y$ and
$\id\colon X\to X$ to denote, respectively, the constant map and the
identity map.  We denote the set of homotopy classes of maps from
$X$ to $Y$ by $[X,Y]$. However, we often do not distinguish between
a map and the homotopy class it represents. In \secref{sec:category}
we frequently consider lifts and sections of certain maps, by which
we mean lifts and sections up to homotopy. If $f \colon A \to B$ is
a map, then $f^*$ denotes pre-composition by $f$ and $f_*$ denotes
post-composition by $f$. Thus we obtain maps of homotopy sets
$f^*\colon [B,X] \to [A,X]$ and $f_*\colon [X,A] \to [X,B]$. We also
use $f^*$ to denote the map induced by $f$ on cohomology and $f_\#$
to denote the map induced on homotopy groups. The \emph{wedge} of
spaces $X$ and $Y$ is denoted $X \vee Y$. The unit interval $[0,1]$
is denoted $I$. The cone on a space $X$ is written $CX$ 
and the path space is $EX$.  The loop space, respectively suspension, of a space
$X$ is denoted $\Omega X$, respectively $\Sigma X$.  We recall that
the standard multiplication $\mu$ on $\Omega Y$ induces a group
structure in $[X, \Omega Y]$ and the standard comultiplication $\nu$
on $\Sigma X$ induces a group structure in $[\Sigma X, Y]$. These
group structures are isomorphic via adjointness, that is, $[X,
\Omega Y] \cong [\Sigma X, Y]$. Under this isomorphism, the identity
$\id \colon \Omega X \to \Omega X$ has an \emph{adjoint} which we
denote $p\colon \Sigma\Omega X \to X$ and the identity $\id \colon
\Sigma X \to \Sigma X$ has an adjoint which we denote by $e\colon X
\to \Omega\Sigma X$. By an \emph{$H$-space}, we mean a pair $(X,m)$
with $m \colon X\times X \to X$ a (not-necessarily associative)
multiplication up to homotopy. Dually, by a \emph{co-$H$-space} we
mean a pair $(X,n)$ with $n \colon X \to X\vee  X$ a
(not-necessarily associative) comultiplication up to homotopy.

\section{Actions on Spaces and Operations on Homotopy Sets}\label{sec:actions}

In this section we consider the \emph{action of one space on
another}, and \emph{cyclic maps from one space to another}.  We
present several general results about these notions that
illustrate the relation between them.  At the end of the section,
we briefly consider the dual notions of coactions and cocyclic
maps.

Let $i_1\colon A \to A\times X$ and $i_2\colon X \to A\times X$
denote the inclusions.

\begin{definition}\label{def:action}
By a (left) \emph{action of $A$ on $X$}, we mean a map $F \colon
A\times X \to X$ that satisfies $F\circ i_2 = \id\colon X \to X$.
If $F\circ i_1 = f\colon A \to X$, then we refer to $F$ as
\emph{an $f$-action}.  Given an $f$-action $F$, we say that $f$ is
\emph{the orbit map of $F$}.  A map $f\colon A \to X$ is called a
\emph{cyclic map} if it is the orbit map of some $f$-action $F$.
We denote the set of cyclic homotopy classes $G(A,X) \subseteq [A,
X]$.
\end{definition}

In the literature, an $f$-action is often referred to as \emph{an
affiliated map (of $f$)}.  There is an alternate point of view on
actions and cyclic maps that is often convenient to adopt: Let $X$
be a space with basepoint $x_0$ and $\map(X,X;\id)$ denote the
component of the identity in the function space of unbased maps from
$X$ to itself. Let $\omega\colon \map(X, X;\id) \rightarrow X$
denote the evaluation map defined by $\omega(h) = h(x_0)$, 
where $x_0$ is the base point of $X$. Then an
$f$-action $F\colon A\times X \to X$ corresponds to a lift
$\tilde{f} \colon A \to \map(X,X;\id)$ of $f$ through $\omega$, that
is, $\omega\circ\tilde{f} = f$.  This correspondence is given by
$\tilde{f}(a)(x) = F(a, x)$. Then the set of cyclic maps $G(A, X)$
may be described as the image of the induced map of homotopy sets
$\omega_*\colon [A, \map(X, X;\id)] \rightarrow [A, X]$.  It follows
easily from this point of view that, if $A$ is a suspension, then
$G(A, X)$ is a \emph{subgroup} of $[A,X]$ \cite{Var}.

We illustrate these notions with some well-known examples.

\begin{examples}\label{exs:actions}
\begin{enumerate}
\item The projection $p_2 \colon A\times X \to X$ is a $\triv$-action.

\item Let $(X, m)$ be an $H$-space.  Then the multiplication
$m\colon X\times X \to X$ is a $\id$-action.

\item Suppose given an $f$-action $F \colon A\times X \to X$ and
any map $h \colon A' \to A$.  Then $F\circ(h\times\id)\colon
A'\times X \to X$ is an $(f\circ h)$-action.  It follows that $h^*
\colon G(A, X) \to G(A', X)$.  By combining this and the previous
observations, we obtain: If $(X, m)$ is an $H$-space, then
$m\circ(f\times\id)\colon A\times X \to X$ is an $f$-action for
any $f\colon A \to X$, and so $G(A, X) = [A, X]$.

\item Let $\omega \colon \map(X,X;\id) \to X$ be the evaluation map
as above.  Define the map $W\colon \map(X,X;\id)\times X \to X$ by
$W(f, x) = f(x)$.  Then $W$ is a (right) $\omega$-action of
$\map(X,X;\id)$ on $X$. Indeed, $W$ is a ``universal action" on
$X$, in that every action of a space $A$ on $X$ factors through
$W$.   From this we see that $\omega$ is a ``universal cyclic
map," in the sense that it is a cyclic map through which every
cyclic map to $X$ factors.

\item Suppose $\xymatrix@1{A \ar[r]^{j} & X \ar[r]^{p} & B}$ is a
fibration.  We call this a \emph{homotopy principal fibration} if
there exists a $j$-action $F \colon A\times X \to X$ of the fibre
on the total space, such that the following diagram homotopy commutes:
$$\xymatrix{A\times X \ar[r]^-{F} \ar[d]_{p_2} & X \ar[d]^{p}\\
X \ar[r]_{p}& B,}$$
see \cite{Pe-Th, Por}. This is the obvious homotopy-theoretic
analogue of the situation in which $A = G$ is a topological group
and $p\colon X \to B$ is a principal $G$-bundle. An important
special case of a homotopy principal fibration is that of a
fibration $\Omega B \to E \to X$ induced from the path-space
fibration over $B$ by some map $g\colon X \to B$.
\end{enumerate}
\end{examples}

An action of one space on another induces an operation on homotopy
sets as follows.  Suppose $F \colon A \times X \to X$ is an $f$-action
and $W$ is any space.  Given $a \in [W,A]$ and $h \in [W,X]$, we
define $a\cdot h \in [W, X]$ as the composition
$$\xymatrix{W \ar[r]^-{\Delta} & W\times W \ar[r]^{a\times h}
& A\times X \ar[r]^-{F} & X,}$$
where $\Delta$ is the diagonal map (see \cite[Ch.15]{Hilton} and
\cite[p.56]{Zab76}). This gives an operation of the set $[W,A]$ on
the set $[W, X]$ with the following properties:

\begin{enumerate}
\item[(i)] $\triv\cdot h = h$;
\item[(ii)] $a\cdot \triv = f_*(a)$;
\item[(iii)] $l^*(a\cdot h) = l^*(a)\cdot l^*(h)$
for any map $l\colon W' \to W$;
\item[(iv)] $F = p_1\cdot p_2$ for the projections
$p_1\colon A\times X \to A$ and $p_2\colon A\times X \to X$.
\end{enumerate}

We now adopt the following notational convention.  If $(X, m)$ is an
$H$-space, we write the binary operation induced on $[W, X]$ by
``$+_m$". On the other hand, if $(W, n)$ is a co-$H$-space, we write
the binary operation induced on $[W, X]$ by ``$\oplus_n$".  This is 
defined as follows: if $a$, $b \in [W,X]$, 
 then $a\oplus _n b = (a|b) \circ n \in [W,X]$.

\begin{example}
Let $(X, m)$ be an $H$-space.   If $f\colon A \to X$ is any map
then, as in (3) of \exsref{exs:actions}, $m\circ(f\times\id)$ is
an $f$-action.  We denote this action by
$$Ac^f_m = m\circ(f\times\id) \colon A\times X \to
X.$$
The induced operation of $[W,A]$ on $[W, X]$ satisfies
$$a\cdot h = f_*(a) +_m h$$
for all $a \in [W, A]$ and $h \in [W, X]$.
\end{example}

We now return to a general $f$-action
$F\colon A\times X \to X$ and identify the operation when $W$ is a co-$H$-space.

\begin{lemma}\label{lem:operation co-H}
Suppose $(W, n)$ is a  co-$H$-space and $F\colon A\times X \to X$
is an $f$-action.  Then we have
$$a\cdot h = f_*(a) \oplus_n h,$$
for all $a \in [W, A]$ and $h \in [W, X]$.
\end{lemma}

\begin{proof}
Observe that the diagram
$$\xymatrix{ & W \vee W \ar[r]^{a\vee h} \ar[dd]^{J} & A \vee X
\ar[dd]_{J} \ar[rd]^{(f\mid\id)} \\
W \ar[ur]^{n} \ar[dr]_{\Delta} & & & X\\
 & W \times W \ar[r]_{a\times h} & A \times X \ar[ru]_{F}}
$$
homotopy commutes, where $J$ is the inclusion.
\end{proof}

\begin{lemma}\label{lem:cyclic on homotopy}
Let $F\colon A\times X \to X$ be an $f$-action and suppose given $a
\in [W, A]$ and $h \in [W, X]$.  Let $(U, n)$ be a co-$H$-space.
Then the induced homomorphism $(a\cdot h)_{*} \colon [U, W] \to [U,
X]$ satisfies
$$(a\cdot h)_{*}(\gamma) = f_{*}a_{*}(\gamma) \oplus_n h_{*}(\gamma),$$
for all $\gamma \in [U, W]$.
\end{lemma}

\begin{proof}
Apply property (iii) above to $(a\cdot h)_*(\gamma) =
\gamma^*(a\cdot h)$, then use \lemref{lem:operation co-H}.
\end{proof}

We may specialize \lemref{lem:cyclic on homotopy} to obtain a
formula for the homomorphism induced on homotopy groups by a map
of the form $a\cdot h$, special cases of which are used in both
\cite{Lu-Op} and \cite{A-L-M}.

\begin{definition}
Let $(X, m)$ be an $H$-space and $F\colon A\times X \to X$ an
$f$-action.  We say that the action $F$ is \emph{$m$-associative}
if the following diagram homotopy commutes:
$$\xymatrix{A\times X\times X \ar[r]^-{F\times \id_X} \ar[d]_{\id_A\times m} & X\times X \ar[d]^{m}\\
A\times X \ar[r]_-{F}& X.}$$
\end{definition}

Notice that in the case in which $A=X$ and the action $F$ is the
multiplication $m$, this definition reduces to that of
homotopy-associativity of $m$.

\begin{proposition}\label{prop:action on an H-space}
Let $(X, m)$ be an $H$-space and $f\colon A \to X$ a map.
\begin{enumerate}
  \item If $m$ is homotopy-associative, then $Ac^f_m$ is an
  $m$-associative $f$-action.
  \item If $F\colon A\times X \to X$ is any $m$-associative
  $f$-action, then $F = Ac^f_m$.
\end{enumerate}
\end{proposition}

\begin{proof}
The first part is clear.  For the second, let $i_1 \colon A \to
A\times X$ denote the inclusion into the first factor and $j\colon
X \to X\times X$ the inclusion into the second factor, so that
$\id_A\times j = i_1\times \id_X \colon A \times X \to A \times X
\times X$. Then we have
\begin{align*}\label{}
  F & = F\circ(\id_A\times m) \circ (\id_A\times j) = m\circ(F\times \id_X) \circ (\id_A\times j) \\
   &  = m\circ\big((F\circ i_1)\times \id_X\big)
  = m\circ(f\times \id_X) = Ac^f_m. \qed
\end{align*}
\renewcommand{\qed}{}\end{proof}

The following result implies that two $f$-actions loop to the same
$\Omega f$-action.  A similar result in a different context appears
in \cite{A-O-S}.  Let $\lambda \colon \Omega A\times \Omega X \to
\Omega(A\times X)$ denote the standard homeomorphism given by
$\lambda(a, b)(t) = \big(a(t), b(t)\big)$, where $t\in I$.

\begin{proposition}\label{prop:loop of an action}
Let $F\colon A\times X \to X$ be an $f$-action and let $\mu$ be the loop
multiplication on $\Omega X$.  Then
$$(\Omega F)\circ\lambda = Ac^{\Omega f}_\mu \colon \Omega A\times \Omega X \to \Omega X,$$
and thus $(\Omega F)\circ\lambda$ is a $\mu$-associative $\Omega
f$-action.  If, in addition,  $ F'\colon A\times X \to X$ is
an $f'$-action such that $\Omega f = \Omega f' \colon
\Omega A \to \Omega X$, then $\Omega F = \Omega F'$.
\end{proposition}

\begin{proof}
Let $p\colon \Sigma\Omega(A\times X) \to A\times X$ be the
canonical map. Then, using (iii) and (iv) of the properties listed above, we have
$$F\circ p = p^*(p_1\cdot p_2) = p^*(p_1)\cdot p^*(p_2),$$
where the
last term uses the operation of $[\Sigma\Omega(A\times X), A]$ on
$[\Sigma\Omega(A\times X), X]$. By \lemref{lem:operation co-H}, it
follows that
$$F\circ p = \Bigl(f\circ p_1\circ p \Bigr)\oplus \Bigl(p_2\circ p\Bigr),$$
where $\oplus$ in
$[\Sigma\Omega(A\times X), X]$ is induced by the suspension structure
on $\Sigma\Omega(A\times X)$.  By taking adjoints and noting that
the adjoint of $p$ is $\id\colon
\Omega(A\times X) \to \Omega(A\times X)$, we have
$$\Omega F = \Bigl((\Omega f)\circ (\Omega p_1) \Bigr)+_\mu \Omega p_2.$$
Let $\pi_1\colon \Omega A\times\Omega X \to \Omega A$ and
$\pi_2\colon \Omega A\times\Omega X \to \Omega X$ be the
projections.  Then
\begin{align*}
(\Omega F)\circ \lambda & =  \Bigl((\Omega f)\circ(\Omega p_1)\circ \lambda \Bigr)+_\mu \Bigl((\Omega p_2)\circ \lambda\Bigr)\\
& =  \Bigl((\Omega f)\circ\pi_1 \Bigr)+_\mu \pi_2 = Ac^{\Omega f}_\mu.
\end{align*}
The other assertions of \propref{prop:loop of an action} follow easily. \qed
\renewcommand{\qed}{}\end{proof}

The following corollary may be known, but we have not found a
proof of it in the literature.  It effectively answers Problem 35
from Stasheff's $H$-space problem list in \cite{H-space}.

\begin{corollary}\label{cor:loop m = loop m'}
If $m, m'$ are two multiplications on $X$, then $(\Omega
m)\circ\lambda = \mu =  (\Omega m')\circ\lambda$.  Consequently,
$\Omega m = \Omega m'$. \qed
\end{corollary}

For the remainder of this section, we consider the dual of the
preceding discussion.  Since most of this can be obtained
\emph{mutatis mutandis} from the previous material, we provide few
details.   Let $p_1\colon X\vee B \to X$ and $p_2\colon X\vee B
\to B$ denote the projections.

\begin{definition}
By a (right) \emph{coaction of $B$ on $X$}, we mean a map $G
\colon X \to X\vee B$ that satisfies $p_1\circ G = \id\colon X \to
X$. If $p_2\circ G = g\colon X \to B$, then we refer to $G$ as
\emph{a $g$-coaction}.  Given a $g$-coaction $G$, we say that $g$
is \emph{the co-orbit map of $G$}.  A map $g\colon X \to B$ is
called a \emph{cocyclic map} if it is the co-orbit map of some
$g$-coaction $G$. We denote the set of cocyclic homotopy classes
by $G'(X,B) \subseteq [X, B]$.
\end{definition}

Notice that there is no natural candidate for the dual of the
function space $\map(X,X;\id)$. So in the coaction setting, we do
not have the alternative point of view that a ``coevaluation map"
would provide.   For results on cocyclic maps, see \cite{Var, Hal-Var,
Lim87, Lee-Woo}. In the case in which $B=K(\pi, n)$, the set of
cocyclic maps $G'(X, B)$ has been called the \emph{$n$th dual
Gottlieb group of $X$} (cf.~\cite{F-L-T}).

We can dualize most of \exsref{exs:actions}; we mention two of
these explicitly.

\begin{examples}\label{exs:coactions}
\begin{enumerate}
\item Suppose given a $g$-coaction $G \colon X \to X\vee B$ and
any map $h \colon B \to B'$.  Then $(\id\vee h)\circ G\colon X \to
X\vee B'$ is an $(h\circ g)$-coaction.  It follows that $h_*
\colon G'(X, B) \to G'(X, B')$.   If $(X, n)$ is a co-$H$-space,
then $(\id\vee g)\circ n\colon X \to X\vee B$ is a $g$-coaction
for any $g$, so $G'(X, B) = [X, B]$.

\item Suppose $\xymatrix@1{A \ar[r]^{j} & X \ar[r]^{p} & B}$ is a
cofibration.  We call this a \emph{homotopy principal cofibration}
if there exists a $p$-coaction $G \colon X \to X\vee B$  such that
the following diagram homotopy commutes:
$$\xymatrix{A \ar[r]^{j} \ar[d]_{j} & X \ar[d]^{i_1}\\
X \ar[r]_-{G}& X\vee B}$$
An important case occurs when the cofibration $A \to X \to \Sigma Z$
is induced from the cofibration $Z \to CZ \to \Sigma Z$ via a map $Z
\to A$ \cite[Ch.15]{Hilton}.
\end{enumerate}
\end{examples}

Suppose $G \colon X \to X\vee B$ is a coaction and $W$ is any
space. Given $a \in [B, W]$ and $h \in [X, W]$, we define $h\cdot
a \in [X, W]$ as the composition
$$\xymatrix{X \ar[r]^-{G} & X\vee B \ar[r]^{h\vee a}
& W\vee W \ar[r]^-{\nabla} & W.}$$
This gives a (right) operation of the set $[B, W]$ on the set $[X,
W]$ \cite[Ch.15]{Hilton}.  We use the same notation for this
operation and the previous one. We note that properties (i)--(iv)
stated earlier may be dualized.

\begin{example}
Let $(X, n)$ be a co-$H$-space.  If $g\colon X \to B$ is any map,
then  $(\id\vee g)\circ n$ is a $g$-coaction which we denote by
$Coac_g^n\colon X \to X\vee B$. In this case, the induced operation
of $[B, W]$ on $[X, W]$ satisfies $h\cdot a = h \oplus_n g^*(a)$,
for all $a \in [B, W]$ and $h \in [X, W]$.
\end{example}

\lemref{lem:operation co-H} dualizes in a straightforward way to
obtain the following.

\begin{lemma}
Suppose $(W, m)$ is an  $H$-space and $G\colon X \to X\vee B$ is a
$g$-coaction.  Then we have
$$h\cdot a = h +_m g^*(a),$$
for all $a \in [B, W]$ and $h \in [X, W]$. \qed
\end{lemma}

\begin{definition}
Let $(X, n)$ be a co-$H$-space and $G\colon X \to X\vee B$ a
$g$-coaction.  We say that $G$ is \emph{$n$-coassociative} if the
following diagram homotopy commutes:
$$\xymatrix{X \ar[r]^-{G} \ar[d]_{n} & X\vee B \ar[d]^{n\vee \id_B}\\
X\vee X \ar[r]_-{\id_X\vee G}& X\vee X\vee B.}$$
\end{definition}

We remark that \propref{prop:action on an H-space} dualizes in a
straightforward way; we omit the dual statement and its proof. Let
$\lambda' \colon \Sigma (X\vee B) \to \Sigma X \vee \Sigma B$
denote the standard homeomorphism.

\begin{proposition}
If $G\colon X \to X\vee B$ is a $g$-coaction and $\nu$ is the
suspension comultiplication on $\Sigma X$, then
$$\lambda'\circ(\Sigma G) = Coac_{\Sigma g}^\nu \colon \Sigma X \to \Sigma X\vee \Sigma B,$$
and thus $\lambda'\circ(\Sigma G)$ is a $\nu$-coassociative $\Sigma
g$-coaction.  In addition, if $G'\colon X \to X\vee B$ is a
$g'$-coaction such that $\Sigma g = \Sigma g'$, then $\Sigma G =
\Sigma G'$. \qed
\end{proposition}

\begin{corollary}
If $n, n'$ are two comultiplications on $X$, then
$\lambda'\circ(\Sigma n) = \nu =  \lambda'\circ(\Sigma n')$, and so
$\Sigma n = \Sigma n'$. \qed
\end{corollary}

We may combine the dual notions discussed above in the following
interesting way.  Suppose given an $f$-action $F\colon A\times X
\to X$ and a $g$-coaction $G \colon X \to X\vee A$ (so $A=B$
here). These give operations of $[X,A]$ and $[A,X]$ on $[X,X]$. In
particular, given any $f \in [A,X]$, $g \in [X,A]$ and
$h \in [X,X]$, we have $h\cdot f,\,  \, g\cdot h
\in [X,X]$ obtained using the coaction and the action,
respectively.

\begin{lemma}\label{lem:cyclic-cocyclic}
With the above notation, $h\cdot f  = g\cdot h $. This
map induces $(g\cdot h)_\# = f_\#g_\# +
h_\#$ on homotopy groups and $(h\cdot f)^* =
h^* + g^*f^*$ on cohomology groups.
\end{lemma}
\begin{proof}  Let $(\id \mid f)\colon X\vee A \to X $ be determined
by $\id $ and $f$, let $(\id, g)\colon X \to X\times A$ be determined
by $\id $ and $g$ and let $T\colon X\times A \to A\times X$ interchange
factors.  Consider the homotopy-commutative diagram
$$\xymatrix{ & X \vee A \ar[r]^{h\vee \id} \ar[dd]^{J} & X \vee A
\ar[dd]_{J} \ar[rd]^{(\id \mid f)} \\
X \ar[ur]^{G} \ar[dr]_{(\id,g)} & & & X.\\
 & X \times A \ar[r]_{h\times \id} & X \times A \ar[ru]_{F \circ T}}
$$
Then
\begin{align*}
h\cdot f & =  ( \id \mid f) \circ (h\vee \id) \circ G  \\ &  =
F\circ T \circ (h\times \id) \circ (\id, g)  =  F\circ (g,h)  =
g\cdot h.
\end{align*}
The other assertions of \lemref{lem:cyclic-cocyclic} follow from \lemref{lem:cyclic on homotopy} and its dual.
\end{proof}

A special case of this situation is
considered in \cite{A-L-M}.

\section{Lusternik-Schnirelmann Category and
Coactions}\label{sec:category}

In this section we relate cocyclic maps to a variant of the
Lusternik-Schnirelmann category of a map. We begin by reviewing the
Ganea fiber-cofiber construction \cite{Gan}.  In this section, we
generally \emph{do} distinguish between a map and its homotopy class.

\begin{definition}\label{def:fiber-cofiber}
Suppose given a fibration $\xymatrix@1{F \ar[r]^{i} & E \ar[r]^{p}
& B}$.  We define a sequence of fibrations
$$\xymatrix@1{F_n(p)
\ar[r]^{i_n} & E_n(p) \ar[r]^-{p_n} & B}$$
inductively for $n \geq 0$. Set $\xymatrix@1{F_0(p) \ar[r]^{i_0} &
E_0(p) \ar[r]^-{p_0} & B}$ equal to the given fibration.  Assume
$\xymatrix@1{F_{n-1}(p) \ar[r]^-{i_{n-1}} & E_{n-1}(p)
\ar[r]^-{p_{n-1}} & B}$ has been defined. Set $E'_{n-1}(p) =
E_{n-1}(p)\cup_{i_{n-1}} CF_{n-1}(p)$, the mapping cone of
$i_{n-1}$.  Define $p'_n \colon E'_n(p) \to B$ by
$p'_n\mid_{E_{n-1}(p)} = p_{n-1}$ and $p'_n\mid_{CF_{n-1}(p)} =
\triv$. Then replace $p'_n$ by a fiber map to obtain a fiber
sequence $\xymatrix@1{F_n(p) \ar[r]^-{i_n} & E_n(p) \ar[r]^-{p_n}
& B}$.
\end{definition}

Note that there is a map $j_{n-1}\colon E_{n-1}(p) \to E_n(p)$
such  that the diagram
$$\xymatrix{E_{n-1}(p)\ar[dr]_{p_{n-1}} \ar[rr]^{j_{n-1}} & & E_n(p)\ar[dl]^{p_n}\\
& B }$$
commutes.

There are two special cases of this construction that we shall
need.

\begin{examples}\label{exs:cofibre-fibre}
\begin{enumerate}
\item Given a space $B$, consider the path space fibration $\Omega B \to EB \to B$.
Applying the above Ganea construction to this fibration, we obtain
a sequence of fibrations which we write as $\xymatrix@1{F_n(B)
\ar[r]^-{i_n} & G_n(B) \ar[r]^-{p_n} & B}$.

\item Given a map $g\colon X \to B$, we first form the fibration
$\Omega B \to E \to X$ induced via $g$ from the path space
fibration over $B$ by the pullback construction.
Then we apply the Ganea construction to this fibration to obtain
a sequence of fibrations which we write as $\xymatrix@1{K_n(g)
\ar[r]^{l_n} & H_n(g) \ar[r]^{r_n} & X}$.
\end{enumerate}
\end{examples}

Notice that the two examples above are related.  From the pullback
construction, there is a map $\tilde{g}\colon E \to EB$ such that
the following diagram commutes
$$\xymatrix{E \ar[r]^{\tilde g} \ar[d]_{p'} & EB \ar[d]^{p}\\
X \ar[r]_{g} & B.}$$
Since the Ganea construction is functorial, we get a map
$\tilde{g_n}\colon H_n(g) \to G_n(B)$ such that the following
diagram is commutative
$$\xymatrix{H_n(g) \ar[r]^{\tilde g_n} \ar[d]_{r_n} & G_n(B) \ar[d]^{p_n}\\
X \ar[r]_{g} & B}$$
for each $n \geq 0$.

The constructions in \exsref{exs:cofibre-fibre} lead to the
following definitions of the category of a map.

\begin{definition}\label{def:cat}
Suppose given a map $g\colon X \to B$.
\begin{enumerate}
\item We say that $\catBG(g) \leq n$ if $g$ can be lifted through $p_n$
to $G_n(B)$.  That is, if there exists $\hat g_n \colon X \to
G_n(B)$ such that
$$\xymatrix{ & G_n(B)\ar[d]^{p_n}\\
X \ar[ur]^{\hat g_n} \ar[r]_{g} & B}$$
commutes up to homotopy.

\item We say that $\catAL(g) \leq n$ if $r_n \colon H_n(g) \to X$ has a
section. That is, if there exists $s_n \colon X \to H_n(g)$ such
that $r_n\circ s_n \simeq \id\colon X \to X$.
\end{enumerate}

We define $\catBG(g)$ (respectively, $\catAL(g)$) to be the least
$n$ such that $\catBG(g)\leq n$ (respectively, $\catAL(g)\leq n$).
\end{definition}

\begin{remarks}\label{rems: catBG}
(1) $\catBG$ of a map is just the Berstein-Ganea category of a map
as defined in \cite{Be-Ga61a}. However, the notion introduced in part (2)
of \defref{def:cat} appears to be new.

(2) We note that $\catBG(\triv) = \catAL(\triv) = 0$, where $\triv : X\to Y$
is the constant map, and
$\catBG(\id_B) = \catAL(\id_B) = \cat(B)$, where $\cat(B)$ denotes
the ordinary Lusternik-Schnirelmann category of the space $B$.

(3)  It follows from the definitions that $\catBG(g) \leq
\catAL(g)$. We will see later that they are generally different.
\end{remarks}

Our main result in this section is that a map $g$ is cocyclic if
and only if it satisfies $\catAL(g) \leq 1$.  To prove this, we need
some notation and a preliminary result.

\begin{notation}
Given a map $g \colon X \to B$, define \emph{the graph of $g$} to be
$\Gamma = \{ (x, g(x)) \mid x \in X \} \subseteq X \times B$. Let
$R_1(g)$ denote the space $E(X\times B; \Gamma, X\vee B)$ of paths
in $X\times B$ that start in $\Gamma$ and end in $X\vee B$ with
fiber map $\tilde\pi_0 \colon R_1(g) \to X$ defined by
$\tilde\pi_0(l) =$ the first coordinate of $l(0)$. One then obtains
$R_1(g)$ as the fiber space over $X$ induced via $(\id,
g)\colon X \to X\times B$ from the fibration over $X\times B$
obtained by converting the inclusion $J \colon X\vee B \to X\times
B$ into a fibration. Finally, if $\gamma$ is a
path in a space $W$ and $s \in I$, then $\gamma_s$ denotes the path
in $W$ defined by $\gamma_s(t) = \gamma(st)$.
\end{notation}

\begin{lemma}\label{lem:tech w}
Suppose given a map $g \colon X \to B$.  Then there is a map
$w\colon H_1(g) \to R_1(g)$ which induces an isomorphism on homology
and is such that the diagram
$$\xymatrix{H_{1}(g)\ar[dr]_{r_{1}} \ar[rr]^{w} & & R_1(g)\ar[dl]^{\tilde\pi_0}\\
& X}$$
commutes.
\end{lemma}

\begin{proof}
The proof, which we sketch here, is an adaptation of an argument
given by Gilbert in [Gil, Prop.3.3]. As stated earlier, $p'\colon E
\to X$ is induced via $g$,
and so we have $\tilde g\colon E \to EB$ (see the discussion after
Examples 3.2). Let $S = \{ (e, \nu) \mid
e\in E, \, \nu\in X^I, \, p'(e) = \nu(0)\}$ and choose a lifting function
$\lambda \colon S \to E^I$ for $p'$.  The extension $p_1'\colon
E\cup C\Omega B \to X$ of $p'$ is given by $p_1'\mid_{E} = p'$ and
$p_1'\mid_{C\Omega B} = \triv$. Then we replace $p_1'$ with a fiber
map $r_1\colon H_1(g) \to X$, where $H_1(g) = \{ (z, \nu) \mid z\in
E\cup C\Omega B,\, \nu\in X^I,\,  p_1'(z) = \nu(1)\}$ and $r_1(z, \nu) =
\nu(0)$.  We now define $w\colon H_1(g) \to R_1(g)$ by
$$w\big( (e, s), \nu\big) = \big( \nu, \big( \tilde g\big( \lambda(e,
-\nu)(1)\big)\big)_s \big),$$
for $s \in I$, $(e, s) \in E\cup C\Omega B \subseteq CE$ and $\nu
\in X^I$ with $p_1'(e, s) = \nu(1)$.  One easily checks that $w$ is
well-defined, has codomain $R_1(g)$, and satisfies $\tilde\pi_0\circ
w = r_1$.  To prove that $w$ is a homology isomorphism, it suffices
to show that $\hat w$, the map induced on fibers, is a homology
isomorphism.  This is done by factoring $\hat w$ through $(\Omega X
\times C\Omega B) \cup (EX \times \Omega B) \subseteq EX\times
C\Omega B$ and showing, as in the argument of Gilbert, that each of
the two resulting maps induces a homology isomorphism.  We omit the
details.
\end{proof}

\begin{theorem}\label{thm:cocyclic vs catAL}
Let $g\colon X \to B$ be a map of simply connected spaces.  Then
$g$ is cocyclic if and only if $\catAL(g) \leq 1$.
\end{theorem}

\begin{proof}
By definition, $g$ is cocyclic if and only if $(\id,g)\colon X \to
X\times B$ factors through $X \vee B$.  One sees easily that this
is so if and only if the fibration $\tilde\pi_0 \colon R_1(g) \to
X$ has a section.  Since $X$ and $B$ are both simply connected,
$w\colon H_1(g) \to R_1(g)$ is a homotopy equivalence by
\lemref{lem:tech w}.  It follows that $g$ is cocyclic if and only
if $r_1\colon H_1(g) \to X$ has a section.
\end{proof}

In the case $X = B$ and $g=\id $, \thmref{thm:cocyclic vs catAL}
reduces to the following well-known result:  $X$ is a co-$H$-space
if and only if $p \colon \Sigma\Omega X \to X$ admits a section
\cite{Gan69}.

The next corollary is a consequence of \remsref{rems: catBG} (3).

\begin{corollary}\label{cor:cocyclic catBG}
If $g\colon X \to B$ is a cocyclic map of simply connected spaces,
then $\catBG(g) \leq 1$.\qed
\end{corollary}

We next give an example to illustrate that the converse of
\corref{cor:cocyclic catBG} does not hold.

\begin{example}
We show that the projection $p_1\colon S^m\times S^n \to S^m$ has
$\catBG(p_1) \leq 1$ and yet is not cocyclic.   That is, we have $1
= \catBG(p_1) < \catAL(p_1)$.  Since $S^m$ is a suspension, there is
a section of $p\colon \Sigma \Omega S^m \to S^m$ and hence a section
of the first Ganea fibration $G_1(S^m) \to S^m$. Thus we have
$\catBG(p_1) \leq 1$. To complete the example, we use the following
fact:  if $g\colon X \to B$ is cocyclic, then for any $\gamma \in
H^{+}(X)$ and $\beta \in H^{+}(B)$, we have $\gamma\cup g^*(\beta) =
0$.  This fact follows easily from the factorization of $(g,\id )$
through the wedge $X \vee B$, where such cup-products are zero
(cf.~the inclusion $\check{G}(X;\Q) \subseteq S(X;\Q)$ of
\cite[Th.1]{F-L-T}). Let $b_m \in H^m(S^m)$ and $b_n \in
H^n(S^n) = H^n(S^m\times S^n)$ be basic classes. Then the cup-product $ b_n\cup
p_1^*(b_m)\not=0 \in H^{m+n}(S^m\times S^n)$, so $p_1$ is not
cocyclic.
\end{example}

We note that there is a third version of the category of a map,
due to Fadell-Husseini \cite{Fa-Hu} and studied by Cornea in
\cite{Cor} as the \emph{relative category} of a map.

\begin{definition}
Suppose given a fibration $\xymatrix@1{F \ar[r]^{i} & E \ar[r]^{p} &
B}$.  Apply the Ganea construction to obtain $\xymatrix@1{F_n(p)
\ar[r]^{i_n} & E_n(p) \ar[r]^{p_n} & B}$ and maps $j_{n-1}\colon
E_{n-1}(p) \to E_n(p)$ such that $p_{n-1} = p_n\circ j_{n-1}$.  We
define $\catFH(p)$ to be  the smallest $n$ such that $p_n$ has a
section $s_n$ with $s_n\circ p = j_{n-1}\circ\cdots\circ j_0$.
\end{definition}

Given a map $g\colon X \to B$, we form the induced fibration
$\xymatrix@1{\Omega B \ar[r] & E \ar[r]^{p'} & X}$ as before.

\begin{lemma}
With the above notation, we have $\catAL(g) \leq \catFH(p')$.
\end{lemma}

\begin{proof}This is immediate from the definition.
\end{proof}

\begin{example}
The inequality above will be strict whenever we find a $p'\colon E
\to X$ which has a section, but which is not a homotopy equivalence.
For instance, let $g = \triv\colon X \to B$.  Then $E = X\times
\Omega B$ and $p'$ is projection onto the first factor.
\end{example}

We have established a connection between a map being cocyclic and
its category. Our point of view is that the category of a map is
inherently of interest, and that this connection motivates interest
in the notion of a cocyclic map. The dual of this connection,
namely, one between the familiar notion of a cyclic map and the
cocategory of a map, then serves to motivate interest in the
cocategory of a map.
\medskip

We summarize some results on the cocategory of a map before
returning to the cocyclic/category side of things.

The fiber-cofiber construction of \defref{def:fiber-cofiber} may be
dualised in the obvious way.  Starting from a cofibration
$\xymatrix@1{X \ar[r]^-{j} & C \ar[r]^-{q} & Z}$, one obtains a
sequence of cofibrations $\xymatrix@1{X \ar[r]^-{j_n} & C_n(j)
\ar[r]^-{q_n} & Z_n(j)}$.  We are interested in two cases of this
construction. If we start from the cofibration $\xymatrix@1{A
\ar[r]^-{j} & CA \ar[r]^-{q} & \Sigma A}$, we obtain a sequence of
cofibrations dual to the Ganea fibrations.  We denote these by
$\xymatrix@1{A \ar[r]^-{j_n} & G_n'(A) \ar[r]^-{q_n} & F_n'(A)}$. If
we start from the cofibration $\xymatrix@1{X \ar[r]^-{j'} & C
\ar[r] & \Sigma A}$, induced via a map $f\colon A \to X$ from
the cofibration $j\colon A \to CA$, then we obtain a sequence of
cofibrations dual to those of (2) of \exsref{exs:cofibre-fibre}.  We
denote these cofibrations by $\xymatrix@1{X \ar[r]^-{j_n} & H_n'(f)
\ar[r]^-{q_n} & K_n'(f)}$. Then we define $\cocatBG(f) \leq n$ to
mean that there exists a map $\tilde f\colon G_n'(A) \to X$ so that
$\tilde f\circ j_n = f$.  We define $\cocatAL(f) \leq n$ to mean
that there exists a retraction of $j_n\colon X \to H_n'(f)$.  From
the definitions, we have $\cocatBG(f) \leq \cocatAL(f)$.

Now it is easy to identify the first Ganea cofibration $A \to
G_1'(A)$ with the canonical map $e\colon A \to \Omega\Sigma A$.
Furthermore, it is well-known that a space $A$ is an $H$-space if
and only if this map admits a retraction.  As we remarked after
\defref{def:action}, a cyclic map factors through the $H$-space
$\map(X,X;\id)$. From this it follows that a cyclic map
$f\colon A \to X$ factors through
$\Omega\Sigma A$, and so $\cocatBG(f) \leq 1$.

On the other hand, we may have a map with $\cocatBG(f) \leq
1$ which is not cyclic.
\begin{example}
Let $f\colon S^3 \to S^3\vee S^3$ be either $i_1$ or $i_2$, the
inclusion into the first or second summand.  Since $S^3$ is an
$H$-space, $f$ factors through $\Omega \Sigma S^3$,
and so $\cocatBG(f) \leq 1$.  But if $f$ is
cyclic, then the Whitehead product $[f, \beta ] =0$, for all
$\beta \in \pi_*(S^3 \vee S^3)$ \cite[Prop.2.3]{Go1}.  Since $[i_1,i_2] \ne 0$,
then $f$ cannot be cyclic.
\end{example}

The main question to be resolved is the following:

\begin{question}
For a map $f\colon A \to X$, is the condition that $f$ be cyclic
equivalent to $\cocatAL(f) \leq 1$?
\end{question}

We now return to cocyclic maps and discuss a result of
Halbhavi-Varadarajan in the light of the above.  In
\cite[Prop.1.2]{Hal-Var}, the following result is shown:

\begin{theorem}\label{thm:H-V}
Let $f\colon A \to X$ be a cyclic map.  Then there is a fibration $X
\to E \to \Sigma A$ such that $f = \partial\circ e$, where
$\partial \colon  \Omega \Sigma A \to X$ is the connecting map of the fibration and $e\colon A \to
\Omega\Sigma A$ is the canonical map.
\end{theorem}

This result makes clear the close relation between cyclic maps and
connecting maps of fibration sequences.  It is an extension of
results of Gottlieb for the case in which $f\colon S^n \to X$ is a
Gottlieb element. The proof makes essential use of the
classifying space for fibrations with given fibre.  As the authors
of \cite{Hal-Var} point out, it is not clear if the dual of this
result holds.  Because there is no analogue of the classifying space
in the setting of cofibrations, their proof cannot be dualized.

\thmref{thm:H-V} implies that every cyclic map $f\colon A \to X$
factors through the canonical map $e\colon A \to \Omega\Sigma A$. We
next give a result which is a slightly weaker version of the dual of
\thmref{thm:H-V}, but which does establish a relation between
cocyclic maps and connecting maps of cofibration sequences.

\begin{theorem}\label{thm:dualH-V}
Let $g\colon X \to B$ be any cocyclic map of simply connected spaces.
Then there is a cofibre sequence
$\Omega B \to E \to H_1(g)$ and a map $s \colon X \to H_1(g)$ such that
$p \circ \partial \circ s = g$,
$$\xymatrix{
H_1(g) \ar[r]^{\partial} & \Sigma \Omega B
\ar[d]^p\\
X \ar[u]^{s} \ar[r]_{g} & B,}$$
where $\partial $ is the connecting map of the cofibre sequence
and $p$ is the canonical map.
\end{theorem}
\begin{proof}
We identify $G_1(B)\equiv \Sigma \Omega B$ with the
mapping cone of $\Omega B \to EB$ and $H_1(g)$ with
the mapping cone of $\Omega B \to E$.  By the discussion
after \exsref{exs:cofibre-fibre}, there is a map $\tilde g_1: H_1(g)
\to \Sigma \Omega B$ such that the following diagram commutes
$$\xymatrix{
H_1(g)\ar[d]^{r_1} \ar[r]^{\tilde g_1} & \Sigma \Omega B
\ar[d]^p\\
X  \ar[r]_{g} & B.}$$
Since $g$ is cocyclic, there exists a section $s\colon X\to H_1(g)$ of
$r_1$ by \thmref{thm:cocyclic vs catAL}.  Thus it suffices to show that $\tilde g_1 = \partial : H_1(g)
\to \Sigma \Omega B$.  But this is a consequence of the third square
in the following mapping of cofibre sequences
$$\xymatrix{\Omega B \ar[d]^{\id} \ar[r] & E \ar[d] \ar[r] &
H_1(g) \ar[d]_{\tilde g_1} \ar[r]^{\partial} & \Sigma \Omega B
\ar[d]^{\id}\\
\Omega B \ar[r] & EB\ar[r] & \Sigma \Omega B \ar[r]_{\id} & \Sigma
\Omega B.&  \qed  }$$

\renewcommand{\qed}{}\end{proof}

We conclude with the observation that the canonical map $p\colon
\Sigma\Omega B \to B$ is universal for cocyclic maps into $B$.

\begin{corollary}\label{cor:cocyclicSigmaOmega}
If $g\colon X \to B$ is a cocyclic map of simply connected spaces,
then $g$ factors through $p\colon \Sigma\Omega B \to B$. \qed
\end{corollary}

\providecommand{\bysame}{\leavevmode\hbox
to3em{\hrulefill}\thinspace}
\providecommand{\MR}{\relax\ifhmode\unskip\space\fi MR }
\providecommand{\MRhref}[2]{%
  \href{http://www.ams.org/mathscinet-getitem?mr=#1}{#2}
} \providecommand{\href}[2]{#2}

\end{document}